# Plato's theory of knowledge of Forms by Division and Collection in the Sophistes is a philosophic analogue of periodic anthyphairesis (and modern continued fractions)

S. Negrepontis

The aim of this paper is to show that Plato's theory of knowledge in the *Sophistes*, obtained by Division and Collection, is a philosophic analogue of the geometric theory of periodic anthyphairesis. Our analysis (in Sections 3 and 4) centers on the two Division-and-Collections—the Being Angler at the start (218b-221c) and the Being Sophist at the end (264b-268d) of the dialogue. Both Beings are shown to be definable and knowable by a process of binary Division, corresponding to the odd steps of a qualitative anthyphairesis, followed by Collection-Logos, a straightforward adaptation of the Logos Criterion for geometric anthyphairesis (described in Section 1). The Logos for the definition of the Sophist is essentially the fundamental analogy of the Divided Line of the *Politeia* 509d-510b. Thus, a Platonic Being (or Form, or Idea) is precisely an entity possessing Division and Collection, namely the philosophic analogue of periodic anthyphairesis; as a result, a Platonic Being is a self-similar One, an entity which is a monad, simple, indivisible and partless in the sense that every part, produced by Division, participates in Logos and is thus, by periodicity, equalized to the whole.

The species, precisely because they are mutually equalized, can serve as units of intelligible numbers, and thus Platonic (eidetic) numbers are generated in every Platonic Being (cf. Sections 5 and 6); indeed, the basic equation described in the *Parmenides* 148d-149d states that the number of different (equalized) units in a Being is equal to (the number of Logoi in a full period)+1, and thus always a finite number. The generation of numbers in the *Parmenides* 143c-144e, the central role of numbers as in the *Philebus* 16c-19b passage, and Aristotle's account of Platonic numbers in the *Metaphysics* (in particular the claim that '[Platonic] numbers are species') that have posed considerable problems of interpretation among modern Platonists, fit perfectly with our anthyphairetic interpretation.

If a Platonic Being is knowable by the periodic anthyphairesis of Division and Logos, then True Judgement turns out to be a finite anthyphairetic approximation of the Platonic Being (Section 7). The descriptions of True Judgement in the *Symposium, Meno, Theaetetus,* and *Sophistes* provide the necessary documentation and support for this anthyphairetic interpretation of True Judgement.

The results of this paper are applied by Negrepontis and Birba-Pappa to the interpretation of Plato's Theory of Falsehood[1].
The interpretation of Pappus' account of Analysis and Synthesis in terms of Plato's dialectics—namely in terms of Division and Collection, and thus, employing the results of the present paper, in terms of periodic anthyphairesis—has been given in Negrepontis-Lamprinidis[2].

This novel interpretation of the *Sophistes* Division and Collection will be expanded in forthcoming publications; we will show (a) that Plato's theory of knowledge is completed in the *Politicus*, where a Platonic Being is shown to possess not simply a periodic, but in fact a palindromically periodic, anthyphairesis, in precisely the same way as the commensurable in power only pairs of line segments (as defined in Book X of Euclid's *Elements*, and corresponding to quadratic irrationals in modern terms), (b) that Theaetetus' mathematical discovery, reported in the *Theaetetus* 147c-148b, is precisely the remarkable theorem of palindromic periodicity for quadratic irrationals, and (c) that the tools for a proof of that theorem are contained in Book X of the *Elements*. Thus, one of the remarkable consequences of the present paper is that Theaetetus must have proven the theorem on the periodicity of the anthyphairetic development of a 'quadratic irrational'.

## 1. PERIODIC ANTHYPHAIRESIS

We outline here the mathematics of 'anthyphairesis', developed by the Pythagoreans, Theodorus, and the geometers, principally Theaetetus, in Plato's Academy and presented, albeit in a highly incomplete manner, in Books VII and X of Euclid's *Elements*.

*(a) Definition of Anthyphairesis.*

> Table 1. Definition of geometric anthyphairesis (Euclid's *Elements* X.2)

Let a, b be two magnitudes (line segments, areas, volumes), with a>b; the anthyphairesis of a to b is the following, infinite or finite, sequence of mutual divisions:

---

[1] S. Negrepontis and S. Birba-Pappa, «Plato's Theory of Falsehood», *this volume*.
[2] S. Negrepontis and S. Lamprinidis, «The Platonic Anthyphairetic Interpretation of Pappus' Account of Analysis and Synthesis», in E. Barbin, N. Stehlikova, C. Tzanakis (eds.), *History and Epistemology in Mathematics Education,* Proceedings of the Fifth European Summer University (2007), Vydavatelsky servis, Plzen, 2008, p. 501-511.

$a = I_0 b + e_1$, with $b > e_1$,
$b = I_1 e_1 + e_2$, with $e_1 > e_2$,
…
$e_{n-1} = I_n e_n + e_{n+1}$, with $e_n > e_{n+1}$,
$e_n = I_{n+1} e_{n+1} + e_{n+2}$, with $e_{n+1} > e_{n+1}$,
…

We set $\mathrm{Anth}(a,b) = [I_0, I_1, \ldots, I_n, I_{n+1}, \ldots]$ for the sequence of successive quotients of the anthyphairesis of a to b.

*(b) Definition (Definitions X.1, 2 of the Elements).* Let a, b be two magnitudes with a>b; we say that a, b are commensurable if there are a magnitude c and numbers n, m, such that a=mc, b=nc, otherwise a, b are incommensurables.
The fundamental dichotomy for anthyphairesis is contained in the following

*(c) Proposition (Propositions X.2, 3 of the Elements).* Let a, b be two magnitudes, with a>b. Then a, b are incommensurable if, and only if, the anthyphairesis of a to b is infinite.

*(d) Anthyphairetic definition of proportion of magnitudes.* Aristotle, in *Topica* 158b-159a, a justly celebrated and extremely important passage for the history of Greek mathematics, refers to a period where no rigorous theory of proportion existed, while in the *Metaphysics* 987b25-988a1, he explicitly states that the Pythagoreans were not conversant with dialectics and 'logoi' (cf. Becker[3]). In the same *Topica* passage, Aristotle tells us that an astounding for its mathematical content (pre-Eudoxian, before Book V of the *Elements)* theory of proportion of magnitudes was discovered, based on the following

*(e) Definition.* Let a, b, c, d be four magnitudes, with a>b, c>d; the analogy a/b=c/d is defined by the condition $\mathrm{Anth}(a,b) = \mathrm{Anth}(c,d)$.
Aristotle in the *Topica* passage gives an important consequence of this anthyphairetic definition of analogy: the following proposition, which Fowler[4] calls:

*(f) The* Topica *Proposition.* If a line segment is divided in a ratio a to b, two orthogonal parallelograms, call them A and B, having base a and b

---

respectively and common height the line segment c perpendicular to the line segments a and b, then $A/B = a/b$.

Proof. Aristotle claims, and it can readily be seen that his claim is true, that obviously

$Anth(A,B) = Anth(a,b)$.

A form of the *Topica* Proposition appears as Proposition VI.1 in Euclid's *Elements*, and is pivotal in applying Eudoxus' theory of analogy, demonstrated in Book V of the *Elements*, to the geometric theory of similarity. The proof of the Proposition in the *Elements* relies on the Eudoxian theory of analogy, but the proof suggested by Aristotle allows, to a considerable extent, the development of much of Book VI of the *Elements*, based on the anthyphairetic definition of analogy.

We will see below that Plato is implicitly employing a philosophic analogue of the *Topica* Proposition—Proposition 4(c1)—to establish the Collection of the Sophist.

*(g) Periodic Anthyphairesis and the Logos Criterion.* An immediate consequence of the anthyphairetic definition of proportion is the following:

Proposition ('the logos criterion' for the periodicity of anthyphairesis). The anthyphairesis of two line segments a, b, with a>b, with notation as in the definition and setting $a=e_{-1}$, $b=e_0$, is eventually periodic, with period from step n to step m-1, if there are indices n, m, with n<m, such that $e_n/e_{n+1} = e_m/e_{m+1}$.

Table 2. Logos Criterion and Periodicity in Geometric Anthyphairesis

Let a, b be two magnitudes (line segments, areas, volumes), with a>b, with anthyphairesis:

$a = I_0 b + e_1$,  with $b > e_1$,

$b = I_1 e_1 + e_2$,  with $e_1 > e_2$,

…

$e_{n-1} = I_n e_n + e_{n+1}$,  with $e_n > e_{n+1}$,

$e_n = I_{n+1} e_{n+1} + e_{n+2}$,  with $e_{n+1} > e_{n+1}$,

…

$e_{m-1} = I_m e_m + e_{m+1}$,  with $e_m > e_{m+1}$,

$e_m = I_{m+1} e_{m+1} + e_{m+2}$,  with $e_{m+1} > e_{m+1}$,

…

so that for some indices n<m we have

$$\boxed{e_n/e_{n+1} = e_m/e_{m+1}}  \qquad \text{(Logos Criterion)}.$$

Then the anthyphairesis of a to b is eventually periodic, and, in fact, Anth(a,b)= $[I_0, I_1, \ldots, \text{period}(I_n, I_{n+1}, \ldots, I_{m-1})]$.

Table 3. Abbreviated Representation
of the Logos Criterion of the Anthyphairesis of a to b

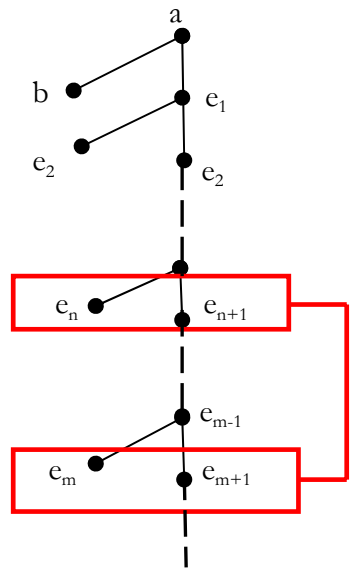

*(h) Reconstruction of the proof of quadratic incommensurabilities by means of the Logos.* There are strong arguments, to appear in detail elsewhere, which follow on from (1) Socrates' exhortation in the *Theaetetus* 145c7-148e5 that Division and Collection imitate the mathematical discovery of Theaetetus relating to incommensurabilities (cf. (i) below), (2) the interpretation of Division and Collection in terms of periodic anthyphairesis in the present paper (Sections 3 and 4, below), and (3) the description, in the *Theaetetus* 147d3-148b2, of the Theodorus-Theaetetus mathematical discoveries in terms of Division and Collection, that the proofs of incommensurabilities given by Theodorus and reported in 147d3-148b2 of square roots of 3, 5,… up to 17, are in fact anthyphairetic, and employ the Logos Criterion (g). Even without any appeal to these powerful correlating arguments, anthyphairetic reconstructions employing the Logos Criterion have been proposed by Zeuthen[5], van der Waerden[6], von Fritz[7], Fowler[8], Kahane[9], a non-

---

[5] H. G. Zeuthen, «Sur la constitution des livres arithmetiques des Elements d' Euclide et leur rapport a la question de l' irrationalite», *Oversigt over det Kgl. Danske Videnskabernes Selskabs Forhandlinger* 5, 1910, p.395-435.

anthyphairetic one by Hardy and Wright[10] and Knorr[11]. In Table 4 below, we outline a reconstruction of the proof of the incommensurability of the line segments a, b, with $a^2=19b^2$, the first one that Theodorus refrains from giving (abbreviated in the sense that the even indexed division-steps are incorporated in the immediately suceeding odd indexed):

Table 4. Abbreviated Anthyphairetic Division and Logos Criterion for $a^2=19b^2$

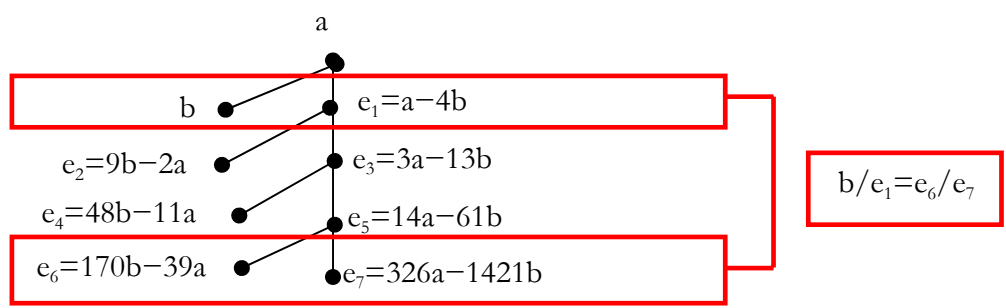

Table 4 is to be understood as follows: we first proceed with the steps of the anthyphairetic Division of a by b, employing elementary computations and expressing at the same time the remainders generated in terms of the initial line segments a and b:

a=4b+$e_1$, with $a_1$<b (hence $e_1$=a-4b),
(and b =2$e_1$+ $e_2$, $e_2$<$e_1$ (hence $e_2$=9b-2a)),
$e_1$ = $e_2$+ $e_3$, $e_3$<$e_2$ (hence $e_3$=3a-13b),
(and $e_2$=3$e_3$+ $e_4$, $e_4$<$e_3$ (hence $e_4$=48b-11a)),
$e_3$= $e_4$+ $e_5$, $e_5$<$e_4$ (hence $e_5$=14a-61b),
(and $e_4$=2$e_5$+ $e_6$, $e_6$<$e_5$ (hence $e_6$=170b-39a)),
$e_5$ =8$e_6$+ $e_7$, $e_7$<$e_6$ (hence $e_7$=326a-1421b); and

we next verify, by direct calculation, the Logos Criterion (indicated in the Table by the coupling of the two expressions in the rectangles), employing the expressions found for the remainders:

b/$e_1$=$e_6$/$e_7$.

It follows that, after the initial ratio a/b, the sequence of successive Logoi

b/$e_1$, $e_1$/$e_2$, $e_2$/$e_3$, $e_3$/$e_4$, $e_4$/$e_5$, $e_5$/$e_6$, $e_6$/$e_7$= b/$e_1$

forms a complete period of Logoi, repeated ad infinitum, and provides full knowledge of the initial ratio a/b, i.e of the quadratic irrational square root of 19, and proving incidentally the incommensurability of the ratio a/b.

---

[6] B. L. van der Waerden, *Science awakening*, translated by A. Dresden, Noordhoff, Groningen, 1954.
[7] K. von Fritz, «The discovery of incommensurability by Hippasus of Metapontum», *Annals of Mathematics* 46, 1945, p. 242-264.
[8] D. Fowler, *The Mathematics…..*
[9] J.-P.Kahane, «La Theorie de Theodore des corps quadratiques reels», *L' Enseignement Mathematique* 31, 1985, p.85-92.
[10] G. H. Hardy and E.M.Wright, *Introduction to the Theory of Numbers*, Oxford, Clarendon Press, 1938.
[11] W. R. Knorr, *The Evolution of Euclidean Elements: A Study of the Theory of Incommensurable Magnitudes and Its Significance for Early Greek Geometry*, Reidel, Dordrecht, 1975.

*(i) Periodic anthyphairesis as a Self-similarity One.* Consider any part, say $e_n$, of a periodic anthyphairesis of a to b. This part $e_n$ then *participates* in the ratio $e_n/e_{n+1}$ (or in the ratio $e_{n-1}/e_n$), and the anthyphairesis of $e_n$ to $e_{n+1}$, being (by periodicity) a cyclic permutation of the anthyphairesis of a to b, essentially coincides with it. In that sense, every part generated in the anthyphairetic process is the same as the whole. Therefore, a pair of magnitudes that possesses periodic anthyphairesis is an example of a self-similar entity. In modern mathematics, such entities exist in abundance, e.g. the Cantor set of the excluded middle third or Sierpinski's gasket, but in ancient Greek mathematics this was the only example with the self-similar property. A self-similar entity clearly deserves the name 'One', since it is everywhere the same. It will become evident in Sections 3, 4, 5, below that Plato considered these periodic anthyphairetic self-similar entities as his model for intelligible Forms and Beings. This is precisely the meaning of Socrates' exhortation in *Theaetetus* 145c7-148e5 to strive to find the true definition of knowledge (of a Platonic Being) by imitating the unifying discovery of Theaetetus and young Socrates on quadratic incommensurabilities (described there in terms of 'powers', 'dunameis').

## 2. THE ANTHYPHAIRETIC INTERPRETATION OF DIVISION AND COLLECTION

*(a) Division and Collection is the philosophic analogue of periodic anthyphairesis.* We will show that periodic anthyphairesis and the Logos Criterion are at the center of Plato's dialectics. The simplest way to see this is to correlate periodic anthyphairesis with Platonic Division and Collection, a method by which Platonic Beings become known to the human soul described in the Platonic dialogues *Parmenides*, *Sophistes*, *Politicus*, *Phaedrus* and *Philebus* (cf. Negrepontis[12]); moreover, the simplest way to grasp the close connection between Division and Collection on the one hand and periodic anthyphairesis on the other is to examine the examples Plato provides of this method in the *Sophistes*. In the present work we will concentrate almost exclusively on the method of

---

[12] S. Negrepontis, «The anthyphairetic nature of Plato's Dialectic», in F. Kalavasis-M. Meimaris (eds.), *Topics in the Didactics of Mathematics V*, Gutenberg, Athens, 2000, p. 15-77 (in Greek); S. Negrepontis, «Plato's theory of Ideas is the philosophic equivalent of the theory of continued fraction expansions of lines commensurable in power only», Manuscript, June 2006; S. Negrepontis, «The periodic anthyphairetic nature of the One in the Second Hypothesis of the Parmenides», Manuscript, September 2005; S. Negrepontis, «The Anthyphairetic Nature of the Platonic Principles of Infinite and Finite» in *Proceedings of the 4th Mediterranean Conference on Mathematics Education*, 28-30 January 2005, Palermo, p. 3-26.

Division and Collection as it is illuminated in the *Sophistes* dialogue, and we will show below (in Sections 3 and 4) that it is a philosophic analogue of the geometric method of establishing periodic anthyphairesis by means of the Logos Criterion, as outlined in section 1.

It follows that a Platonic Form is One and Many, not in the sensible cumulative sense of One man with many limbs (specifically rejected by Plato in the *Parmenides* 128e-130a and in the *Philebus* 14d-e), but in the intelligible sense that the Platonic Being is Many, and in fact infinitely Many, in the self-similar sense that it has infinitely many parts, but is still essentially a partless One, in the sense that each part-species is the same as the whole. The self-similarity of Platonic Beings is examined in Section 5.

The equalization of the parts-species effected by periodicity and self-similarity opens the way for an understanding of Platonic ('eidetic') numbers, in that it produces a finite number (in fact, a number equal to the length of the period plus one) of equalized units, which are precisely the parts-species within a period. Platonic numbers are examined in Section 6.

The anthyphairetic interpretation of Division and Collection opens the way for an understanding of True Judgement, since a Platonic Being is also knowable as True Judgement plus Logos. True Judgement, as examined in Section 7, turns out to be the philosophic analogue of any finite anthyphairetic approximation of the Platonic Being.

*(b) Earlier efforts to connect Plato's philosophy with anthyphairesis.* A few students of Plato have detected in some parts of Plato's writings some connection with the arithmetic/geometric concept of anthyphairesis. To my knowledge, they are the following:

(i) Alfred E. Taylor[13] and (ii) D'Arcy W. Thompson[14],
who saw a connection between the (anthyphairetic) side and diameter numbers and the Excess and Defect described in the *Epinomis*; no wider implications were, however, realized;
(iii) Charles Mugler[15],
who sensed the connection between geometric anthyphairesis via the *Theaetetus* 147-8 mathematical passage, and the *Sophistes*' Division and Collection, but in a defective way; Cherniss[16], who in his 1951 book review 'demolished' Mugler's approach, but for the wrong reasons; and,
(iv) Jules Vuillemin[17],
who sensed correctly that Platonic Division and Collection is related to periodic anthyphairesis, though the connection he envisaged was not the correct one. In addition,
(v) David Fowler[18],
who suggests that anthyphairesis was important in Plato's Academy, but does not explain in which way.

---

[13] A. E. Taylor, «Forms and Numbers: A Study in Platonic Metaphysics», *Mind* 35, 1926, p.419-440; ibid., 36, 1927, p.12-33.

[14] D' Arcy W. Thompson, «Excess and defect: Or the little more and the little less», *Mind*, 38, 1929, p.43-55.

[15] C. Mugler, *Platon et la Recherche Mathematique de son Epoque*, Editions P. H. Heitz, Strasbourg-Zurich, 1948.

[16] H. Cherniss, «Plato as Mathematician», *Review of Metaphysics*, 4, 1951, p. 395-425.

But these views remained decidedly marginal; researchers like Mugler and Vuillemin were unable to convince Platonists of the importance of anthyphairesis in Plato's work, while, conversely, Platonists were—and still are—unable to grasp that the Platonic method of Division and Connection describes Platonic Being. In my opinion, this double failure is due to the failure to understand the manner in which Collection turns the infinitely many parts of the anthyphairetic Division into an entity that deserves to be called One, and thus be conceived as a Platonic Being. Precisely, this will be the main task of the present paper, realized in Sections 3, 4, and 5, below.

## 3. THE DIVISION AND COLLECTION OF THE PLATONIC BEING 'THE ANGLER' IN THE *SOPHISTES* 218b-221c

*(a) The Division of the 'Angler'.* The method of Division and Collection, also called 'Name and Logos' (cf. *Theaetetus* 201e2-202b5, *Sophistes* 218c1-5, 221a7-b2, 268c5-d5), is exemplified at the start of the *Sophistes* 218b-221c by the definition of the Angler. In the scheme below, we reproduce the binary division process which leads to the Angler.

---

[17] J. Vuillemin, *Mathematiques Pythagoriciennes et Platoniciennes*, Albert Blanchard, Paris, 2001.
[18] D. Fowler, *The Mathematics…..*

Table 5. The Division of the Angler (*Sophistes* 218b-221c)

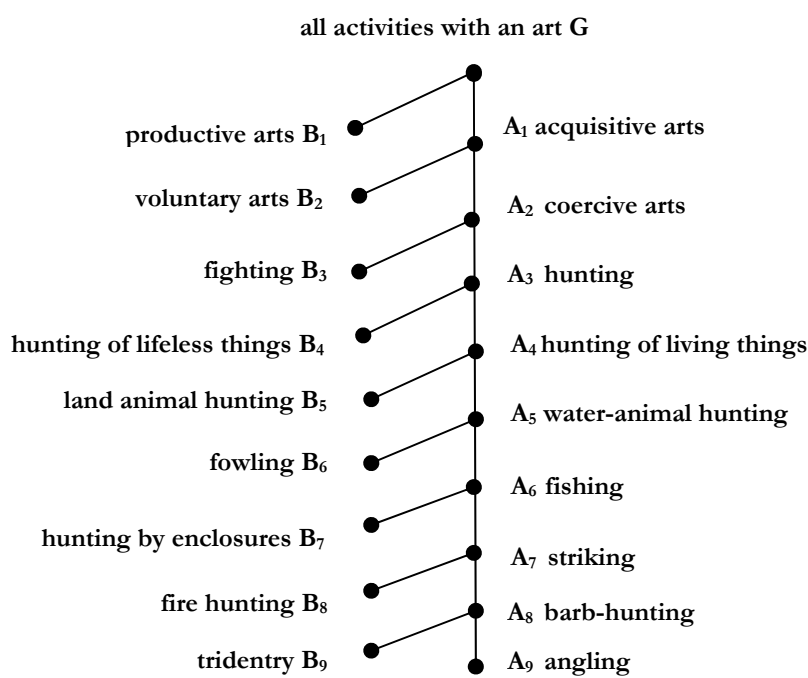

*(b) Collection-Logos of the Angler.* The description of the Logos-Collection of the Angler is contained in the *Sophistes* 220e2-221c3 passage (which we have divided into two parts [A] and [B] for the sake of convenience):

> [A] 'Stranger: Then of striking which belongs to barb-hunting, that part which proceeds
> downward from above ('anothen eis to kato'),
> is called, because tridents are chiefly used in it, tridentry, I suppose….
> Stranger: The kind that is characterized by the opposite sort of blow, which is practised with a hook and strikes,... and proceeds
> from below upwards ('katothen eis tounantion ano'),
> being pulled up by twigs and rods.
> By what name, Theaetetus, shall we say this ought to be called?
> Theaetetus: I think our search is now ended and we have found the very thing we set before us a while ago as necessary to find.
> Stranger: Now, then, you and I are
> only agreed about the name of angling,' (220e2-221b1)

> [B] 'but we have acquired also a satisfactory 'Logos' of the thing itself. For ('gar')
> of art as a whole, half was acquisitive,
> and of the acquisitive, half was coercive,
> and of the coercive, half was hunting,
> and of hunting, half was animal hunting,
> and of animal hunting, half was water hunting,
> and of water hunting

> the whole part from below ('to katothen tmema holon')
> was fishing,
> and of fishing, half was striking,
> and of striking, half was barb-hunting,
> and of this the part in which the blow is pulled
> from below upwards ('peri ten katothen ano') was angling.'[19] (221b1-c3).

In [A], the opposing relation of Tridentry to Angling is carefully explained: all Fishing with a hook is divided into
Tridentry (=Fishing with a trident), which is described as
Fishing with a hook with an art that proceeds from above downwards, and
Angling (=Fishing with a rod), which is described as
Fishing with a hook with an art that proceeds from below upwards.
We have pushed Division all the way to the Angling; thus, we have certainly found 'the name' of Angling. But now, in [B], it is claimed that 'the Logos' of the Angling has been found, too.
The justification—the proof—that we have indeed found the Logos, too, is contained in the remainder of [B], since this remaining part of [B] starts with a 'for' ('gar'), and this justification can be seen to consist of:
(i) an accurate recounting of all the division-steps, abbreviated in the sense that of the two species into which each genus is divided, only the one that contains the Angler is mentioned, while the species opposite to it is omitted;
(ii) a reminder that the last species, the angling, is characterised as the part of its genus that proceeds 'from below upwards'; and,
(iii) the ONLY new information (since (i) and (ii) are repetitions of things already contained in the Division and in [A]), which concerns the species of fishing, three steps before angling, and which informs us for the first time that this species is 'the whole part from below' of its genus.
Since this is an abbreviated account, there is no explicit information on the species opposite to 'fishing'—namely 'fowling', but since 'fishing' was described not simply as 'the part from below' of its genus, but emphatically as 'the whole part from below', it follows that the opposite species—'fowling'—must be characterised as 'the (whole) part from above' of the same genus. In fact, there can be no other justification for the presence of the term 'whole' in the description of 'fishing' with a view to arguing that we have obtained 'Logos', except to indicate and imply this description for its opposite species, 'fowling'.
We recall that the part of [B] from the word 'for' ('gar') is explicitly a justification of the claim that we have succeeded in finding the 'Logos' of the Angling. We may then ask: what is the 'Logos' of the Angler that reasonably results from such a justification? There can really be only one answer: the 'Logos' we are looking for is the equality of the 'philosophic ratio' of Tridentry to Angling, namely the equality of the ratio 'of from above downwards to from below upwards', to the ratio 'of Fowling to Fishing'.

---

[19] Based on Plato, *Theaetetus and Sophist*, translated by H.N. Fowler, Loeb Classical Library, Cambridge, Mass., 1921.

Since the species Tridentry and Angling forms a pair of opposite species, and the species Fowling and Fishing form another pair of opposite species in the Division Scheme for the Angler, the resulting 'Logos' bears a most uncanny similarity to the Logos Criterion for the periodicity of the anthyphairesis of geometric magnitudes, and of geometric powers in particular.

*(c) The Division and Collection of the Angler.* The Division and Collection of the Angler thus takes the following form:

Table 6. Division and Collection of the Angler (*Sophistes* 218b-221c)

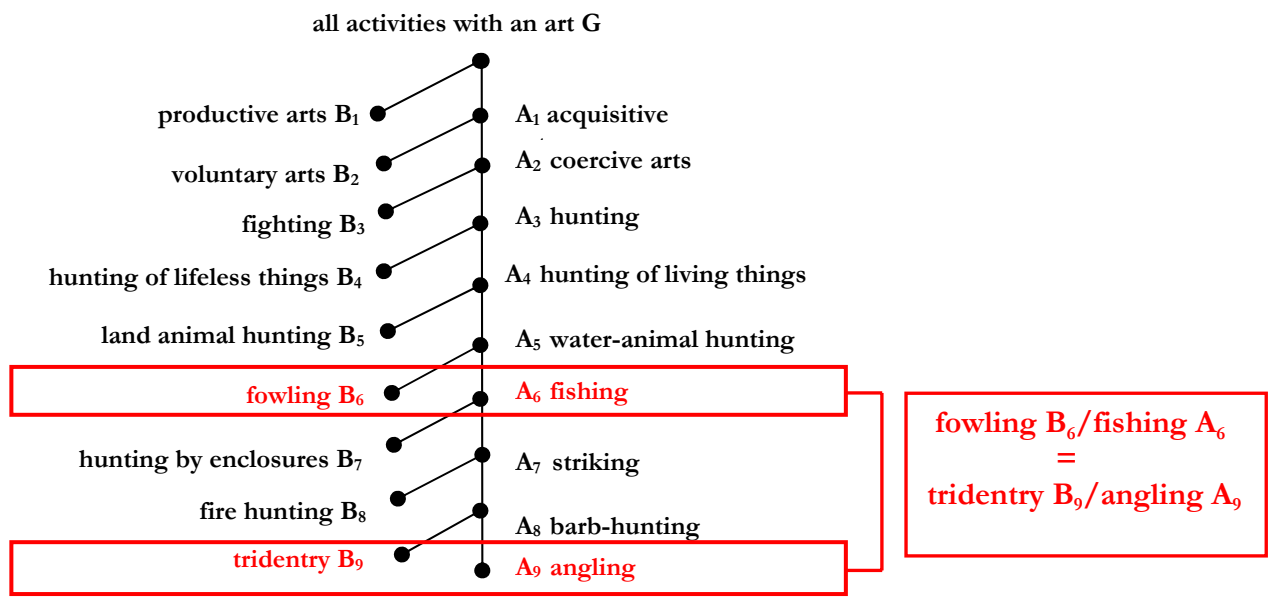

Thus the Division and Collection of the Angler consists of the Division, described in a, which is analogous to the abbreviated anthyphairetic model, as given in 1a, and of the Logos, described in b, analogous to the Logos Criterion for periodicity for geometric anthyphairesis examined in 1(g),(h).

## 4. THE DIVISION AND COLLECTION OF THE PLATONIC BEING 'THE SOPHIST' IN THE *SOPHISTES* 234e-236d & 264b-268d

It now appears that the Division and Collection of a Platonic Being—and the Angler is certainly a lowly paradigm of a Platonic Being—is very much like the anthyphairetic Division and the Logos Criterion of a geometric 'power'. When Socrates expressed, in the *Theaetetus* 145c7-148e5, his exhortation to imitate the geometric situation, it seems that he meant a much closer imitation that anybody has thus far suspected! But before proceeding on to wide-ranging conclusions, it would be prudent to examine whether the Division and Collection of the Sophist, in the *Sophistes* 264b-268d, has the same kind of structure and, in particular, whether there is a similar type of 'Logos'. So we shall take a look at the Division and Collection of the Sophist.

*(a) The Division of the Sophist.* The Division for the Sophist follows the same pattern as the Division of the Angler, starting with a Genus, in this case 'all the productive arts', proceeding by binary division of each Genus to two species, where the next Genus is that species of the previous step in the Division that contains the entity to be defined, in this case the Sophist, and ending with the division-step that produces the Sophist as a species. The whole division scheme is as follows:

Table 7. Division of the Sophist (*Sophistes* 264b-268d)

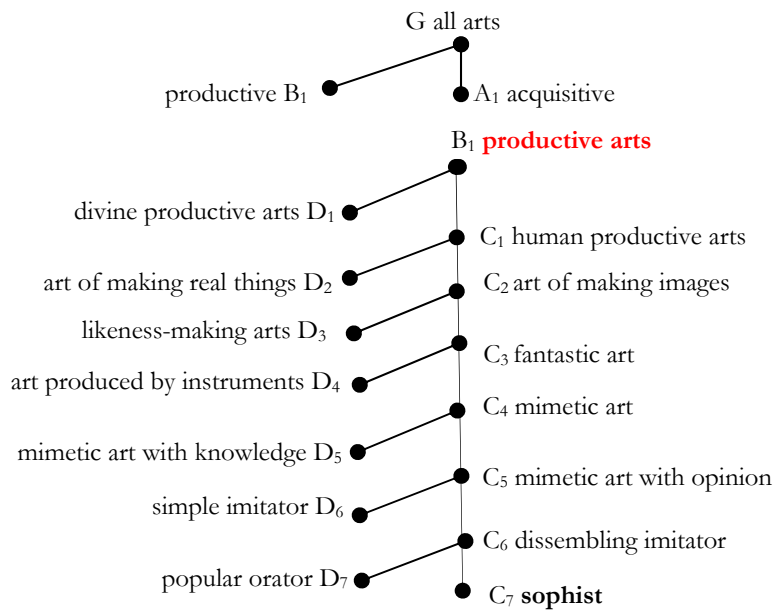

*(b) The fundamental analogy in the Divided Line of the Politeia 509d-510b.* We will now prepare the ground for the Logos-Collection of the Sophist. The fundamental analogy in the Divided Line of the *Politeia* 509d-510b plays a central role in the Logos Criterion of the Sophist. Here is the passage:

| 'Conceive then, said I, as we were saying, that there were two entities, and that | |
|---|---|
| | one of them is sovereign over the intelligible order and region |
| and the other over | |

| | | |
|---|---|---|
| | the world of the eyeball | |
| | …You surely apprehend the two types, | |
| the visible ('horaton') | | |
| | | and the intelligible ('noeton') |
| …Represent them then, as it were, by a line divided into two unequal sections and cut each section again in the same ratio ('ana ton auton logon') | | |
| (the section, that is, of the visible | | |
| | | and that of the intelligible order), |
| and then as an expression of the ratio of their comparative | | |
| | | clearness |
| and obscurity | | |
| | you will have, as | |
| one of the sections of the visible world, images ('eikones'). By images I mean, first, shadows, and their reflections in water and on surfaces of dense, smooth and bright texture, and everything of that kind… | | |
| | As the second section assume that of which this is a likeness or an image ('ho touto eoike'), that is, the animals about us and all plants and the whole class of objects made by man… | |
| | Would you be willing to say, said I, that the division in respect of | |
| | | truth |
| or the opposite | | |
| | | ('aletheia te |
| kai me') | | |

|  | is expressed by the proportion: |  |
|---|---|---|
|  | as is |  |
| the opinable ('doxaston') |  |  |
|  |  | to the knowable ('gnoston') |
|  | so is |  |
| the likeness ('to homoiothen') |  |  |
|  | to that of which it is a likeness ('to ho omoiothe')? |  |
|  | I certainly would.'[20] |  |

This analogy in the *Politeia* Divided Line is rendered as follows:

| Let L be a line (segment), and divide the line L into two unequal sections, say | |
|---|---|
| | A |
| and B, | |

with

|  | A representing the intelligible domain, |
|---|---|
| and B the visible, or rather the sensible, domain. | |

(Incidentally, the construction of the division of a line segment into a given ratio is contained in Proposition VI.10 of the *Elements*).

|  | Then divide section A into two sections, say | |
|---|---|---|
|  |  | C |
|  | and D, |  |
| and divide section B into two sections, say |  |  |
|  | E |  |
| and F, |  |  |
| in such a way that B/A=D/C=F/E. Further, | | |
| F represents the images in the sensible domain B, |  |  |
|  | and E |  |

---
[20] Based on Plato, *The Republic*, translated by P. Shorey, Loeb Classical Library, Cambridge, Mass.,1935.

|  | the entities |  |
|---|---|---|
| in the sensible domain B to which |  |  |
| these are images, |  |  |
|  | i.e. the real entities |  |
| in the sensible domain B. |  |  |
|  |  | Also the intelligible domain A is identified with the domain of knowledge, |
| and the sensible domain B is identified with the domain of opinion. |  |  |
| Hence the following proportion holds: The ratio of the opinable B to the knowable A is equal to the ratio of the images F to the real entities E. |  |  |

This, then, is the fundamental analogy of the Divided Line in the *Politeia* 509d-510b:

the opinable/the knowable = the likeness/that of which it is a likeness.

Equivalently: real things/ images = knowable/opinable.

*(c) Rendering the fundamental analogy of the Divided Line usable for the Division and Collection of the Sophist (Sophistes 265e8-266d7).* More than one fifth (the part 265e8-266d7) of the whole description of the Division and Collection of the Sophist (264b-268d) is devoted to considerations that would appear unnecessary and superfluous at first sight.

In fact, the division for the Sophist starts with the Genus of all productive arts, all arts that produce something; this genus is divided into two species, the divine producing arts and the human producing arts. Since Sophistry is one of the human producing arts, the species human producing arts becomes the next Genus to be divided. This step is justified, considered quite reasonable, and finally accepted at 265e7.

The second division-step divides the human producing arts into human arts producing real things and human arts producing images. The full justification offered for this step is included in the passage 266c7-d4, where the construction of a house is cited as an example of a human art producing a real thing, the painting of a house as an example of a human art producing images, and the division into the two species is stated.

But Plato includes an additional argument, superfluous to the Division,

which is included in the following passage 265e8-266b1 (the previous two sentences (265e3-7) are included for purposes of continuity):

> 'But I will assume that things which people call natural are made by divine art,
> and things put together by man out of those as materials are made by human art,
> and that there are accordingly two kinds of art
> the one human and
> the other divine.
> Theaetetus: Quite right.' (265e3-7).

> 'Stranger: Now that there are two,
> divide each of them again.
> Theaetetus: How?
> Stranger: You divided all productive art widthwise, as it were, before;
> now divide it lengthwise.
> Theaetetus: Assume that it is done.
> Stranger: In that way we now get four parts in all;
> two belong to us and are human,
> and two belong to the gods and are divine.
> Theaetetus: Yes.
> Stranger: And again, when the section is made the other way,
> one part of each half has to do with the making of real things, and
> the two remaining parts may very well be called image-making;
> and so productive art is again divided into two parts.'[19] (265e8-266b1).

A precise description of the geometrical construction of this passage is as follows:

We represent all productive arts by a square, say P, of side p, with vertices K,L,M,N (so each of the sides KL, LM, MN, NK is equal to p) (see Table 8 below).

We divide P widthwise by bringing the horizontal line RS parallel to the base KL, so that all productive art is divided into two rectangles, one KLSR representing human productive arts—call it H, and the other RSMN representing divine productive arts—call it D.

We divide P lengthwise by bringing the vertical line TU parallel to the side KN, so that all productive art is divided into two rectangles, one KTUN representing all productive arts producing real things—call it R, and the other
TLMU representing all productive arts producing images—call it l.

The initial genus of all productive arts is represented not as a line segment (as would be natural), but uniquely and unexpectedly as a square; it is clear that it is thus represented so as to be divisible simultaneously in two independent ways, horizontally and vertically: horizontally in terms of the first division (divine/human), vertically in terms of the second division (real things/images).

This representation serves no purpose whatsoever for the Division of the Sophist, where the species opposite the species to which the Sophist belongs always stay undivided, and does not in fact influence the subsequent division-steps.

Since the Division could serve no possible purpose, we conclude that the entire argument must have some role to play in the Collection-Logos. When we think in terms of Collection and Logos, this whole argument, which would otherwise appear out of place, immediately assumes a rationale:

in fact the ratio of the second step in the Division of the Sophist is:

the ratio of human productive arts producing <u>real entities</u> to human productive arts producing <u>images</u>.

How can this ratio be exploited for the Collection-Logos Criterion? At this point, we can hardly fail to recall the quite similar ratio which appears in the analogy of the Divided Line and was examined in detail in b:

'The ratio of the opinionable to the knowable is equal to the ratio of <u>images</u> to <u>real entities</u>'.

Inverting the ratios, we have:

the ratio of the knowable to the opinionable is equal to the ratio of real entities to images.

This dialectical Logos would certainly serve excellently as the Logos in a Division and Collection, while the second ratio in the Division of the Sophist is beguilingly close to the ratio of the Divided Line, albeit with a rather sophisticated difference: the second ratio in the definition of the Sophist is a ratio which concerns only the human productive arts, while the Divided Line ratio is a ratio which concerns all the productive arts, both divine and human. This could be put down to carelessness, but since Plato is intent on adhering to the highest standards of precision, he must prove something if he is to rely on it. Thus, the passage 265e8-266b1, which seemed unnecessary and out of place, is exactly what is needed to prove the following:

*(c1) Proposition.* The ratio of the second step in the Division of the Sophist, namely the ratio of the human productive arts producing real things to the human productive arts producing images is equal to the ratio of all productive arts producing real things to all productive arts producing images in the *Politeia* Divided Line

Proof. We employ the above notation (see also Table 8).
Of course, the horizontal line RS cuts the side KN into two segments, KR=h and RN=d, and the vertical line TU cuts the side KL into two segments, KT=r and TL=i.
If V is the point where the horizontal line RS meets the vertical line TU, then
all productive arts producing real things=R=the rectangle KTUN with sides r and p= rp, and
all productive arts producing images=I=rectangle TLMU with sides i and p= ip; and
human productive arts producing real things= rectangle TVRK with sides r and h= hr, and
human productive arts producing images=rectangle TLSV with sides I and h= hi.
Then
the ratio (1) of the human productive arts producing real things to the human productive arts producing images
= rectangle TVRK / rectangle TLSV
= hr/hi,
and
the ratio (2) of all productive arts producing real things to all productive arts producing images
= rectangle KTUN/ rectangle TLMU
=R/I=rp/ip.

Now we employ the fundamental *Topica* Proposition 1(f), essentially Proposition VI.1 of the *Elements*, according to which
the ratio (1)= hr/hi= r/i= rp/ip =R/I = the ratio (2),
so that the Proposition is proved.

Table 8.

**r/i**=arts producing real things/arts producing images =
human arts producing real things/human arts producing images=**hr/hi**

| N                                      | U                              M |
|----------------------------------------|----------------------------------|
| **d divine** arts producing **real** things<br>**dr** | **divine** arts producing **images**<br>**di** |
| R                                      | V                              S |
| **h human** arts producing **real** things<br>**hr** | **human** arts producing **images**<br>**hi** |
| K                **r**                 | T                **i**         L |

*(d) The Collection-Logos of the Sophist.* The Logos Criterion of the Sophist does not appear obvious or simple-minded as with the Angler. But the Angler, after all, was chosen as an example not because we had some great philosophical interest in the Art of Angling, but precisely because there was a simplicity, even a naivety, to the Division, but more so in the Logos. From the Sophist, we expect something less obvious.

Once we notice the proximity of the ratio in the second division-step to one of the ratios (images to real things) in the analogy of the Divided Line, and that the argument examined in (c) seems to have been inserted by Plato not only for the purpose of precision but also—and still more so—for the purpose of bringing attention to this ratio, we cannot fail to note that the ratio of the fifth step in the Division of the Sophist has the same type of proximity to the other ratio (the opinionable to the knowable) in the analogy of the Divided Line. The fifth division of the Genus, mimetic / imitative art, into the two species of imitative art with knowledge and imitative art with opinion is contained in the *Sophistes* 267a10-e3 passage below (which also occupies more than a fifth of the entire description of the Division and Collection of the Sophist):

> 'Stranger: But it is surely worth while to consider, Theaetetus, that the mimetic art ('mimetikon') also has two parts… Some who imitate do so with knowledge ('oi men eidotes') of that which they imitate, and others without such knowledge ('oi de ouk eidotes'). And yet what division can we imagine more complete than that which separates ignorance and knowledge ('agnosias kai gnoseos')?
> Theaetetus: None.
> …
> Stranger:… Are there not many who have no knowledge of it, but only a sort of opinion ('agnoountes men, doxazontes de')…
> Theaetetus: Yes, there are very many such people.
> …
> Stranger: Then I think we must say that such an imitator ('mimeten') is quite distinct (heteron) from the other, the one who does not know from the one who knows ('ton agnoounta tou gignoskontos').
> Theaetetus: Yes.
> …
> Stranger:…call the imitation which is based on opinion, opinion-imitation ('ten men meta doxes mimesin doxomimetiken'), and that which is founded on knowledge, a sort of scientific imitation ('ten de met' epistemes historiken tina mimesin').'[19]

It is clear that the same method employed for the ratio of the second division-step can be used to prove the following Proposition:

*(d1) Proposition.* The ratio of the fifth step in the Division of the Sophist—namely, the ratio of the imitative arts with knowledge to the imitative arts with opinion—is equal to the ratio of knowledge to opinion of the *Politeia* Divided Line.

Plato does not make any explicit reference to the procedure, but seems to implicitly ascertain the equality of these two ratios in 267b7-9, where the relation of imitation to knowledge and imitation to ignorance is equated to the relation of knowledge itself to ignorance itself.

We are now ready to describe the Logos Criterion for the Sophist:

*(d2) Proposition (Logos of the Sophist).* The ratio of the second step in the Division of the Sophist, namely the ratio of the human productive arts producing real things to the human productive arts producing images, is equal to the ratio of the fifth step in the Division of the Sophist, namely the ratio of the imitative arts with knowledge to the imitative arts with opinion.

Proof. We use the analogy of the Divided Line, explained in b, and Propositions (c1) and (d1).

*(e) The Division and Collection of the Sophist.* The complete Division and Collection of the Sophist can thus be summarised in the following Table 9 (in which the other two connections shown, to be ignored at the moment, will be explained in (f), below):

Table 9. Division and Collection of the Sophist (*Sophistes* 264b-268d)

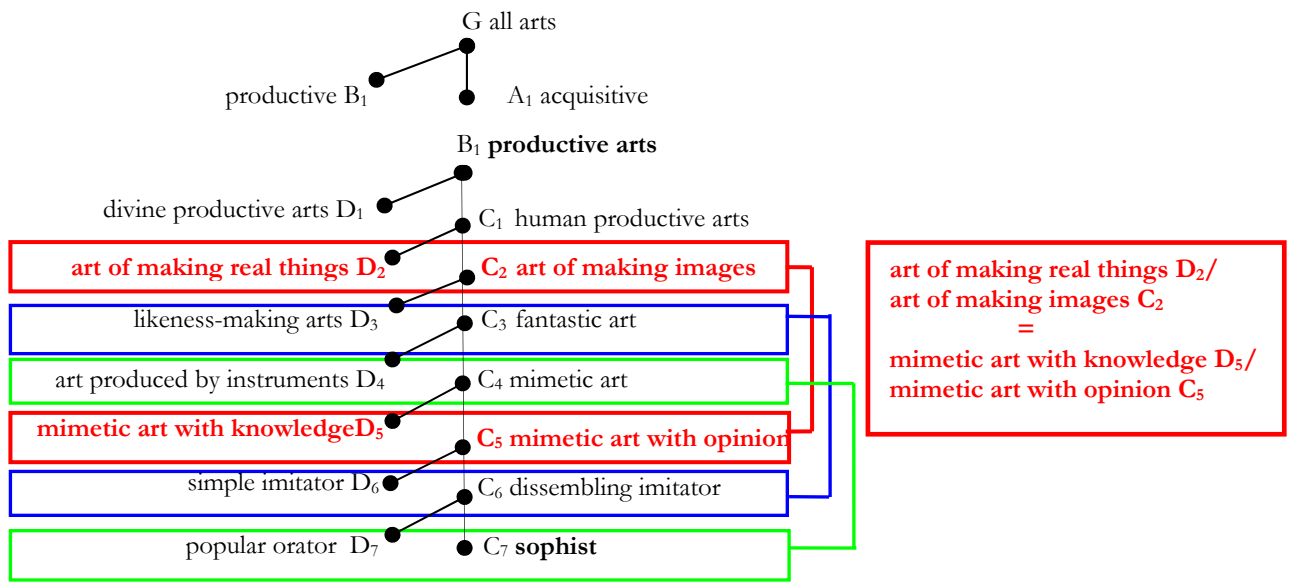

Thus, we have full confirmation for the anthyphairetic meaning of Logos that we proposed in the Division and Collection of the Angler. This time, the Logos is not of the simple-minded, obvious type that appeared in the definition of the Angler (from above downwards/from below upwards), but of a more sophisticated and philosophical kind, playing a central role in the dialectics of the *Politeia*[21].

*(f) Philosophical consequences for the next two Logoi of the Division and Collection of the Sophist analogous to the geometrical consequences of periodic anthyphairesis.* The definition of the Sophist presents an additional feature that ties it even closer to the mathematical model: there are two Logoi after the Logos Criterion, and if the mathematical

---
[21] The correlation of the *Politeia* ratio in the divided line with the Division and Collection of the Sophist in the *Sophistes* is hinted at by Proclus in *eis Politeian* 1,290,7-10.

anthyphairetic model is indeed followed, we would expect to have two further equalities of Logoi: namely the logos of the third step should be equal to the ratio of the sixth division-step, and the logos of the fourth division-step should be equal to the ratio of the seventh and final division-step**.**

*(f1) The equality of the third and sixth ratios of the Division of the Sophist.* In the third division-step, the Genus 'humanly produced images' is divided into two species,

> 'Stranger: I see the likeness-making art ('eikastiken') as one part of [the image-making art]. This is met with, as a rule, whenever anyone produces the imitation by following the proportions of the original in length, breadth, and depth, and giving, besides, the appropriate colors to each part.' (235d6-e2);

> 'Stranger: Now then, what shall we call that which appears, because it is seen from an unfavorable position, to be like the beautiful, but which would not even be likely to resemble that which it claims to be like, if a person were able to see such large works adequately? Shall we not call it, since it appears, but is not like, an appearance?
> Theaetetus: Certainly.' (236b);

> 'Stranger: And to the art which produces appearance, but not likeness ('phantasma all' ouk eikona'), the most correct name we could give would be "fantastic art," ('phantastiken') would it not?
> Theaetetus: By all means.
> Stranger: These, then, are the two forms of the image-making art ('eidolopoiikes') that I meant, the likeness-making ('eikastikes') and the fantastic ('phantastikes').' (236c3-7);

> 'Stranger: We must remember that there were to be two parts of the image-making class ('eidolourgikes'), the likeness-making ('eikastikon') and the fantastic ('phantastikon')'[19] (266d8-9).

In the sixth division-step, described in 267e7-268a8, the Genus opinionated-imitator is divided into the two species, 'the simple' ('haploun') and 'the dissembling' ('eironikon') opinionated-imitator, the simple imitators being those who
'think they know that about which they have only opinion',
while the 'dissembing' are those who
'because of their experience in the rough and tumble of arguments, strongly suspect and fear that they are ignorant of the things which they pretend before the public to know'.

Thus the simple imitators do not distort their opinion, but rather express a likeness of their opinion, while the dissemblers distort and disguise their opinion behind a false appearance. Therefore we have

*(f2) Proposition.* The ratio of the third step in the Division of the Sophist, namely the ratio of the likeness-making arts to the phantastic arts, is equal to the ratio of the sixth step in the Division of the Sophist, namely the ratio of the simple imitator to the dissembling imitator.

*(f3) The equality of the fourth and seventh ratios.* In the fourth division-step, described in 267a1-b3, the Genus, fantastic arts, is divided into two species, as follows:

> 'Stranger: Let us, then, again bisect the fantastic art.
> Theaetetus: How?
> Stranger: One kind is that produced by instruments, the other that in which the producer of the appearance offers himself as the instrument.
> Theaetetus: What do you mean?
> Stranger: When anyone, by employing his own person as his instrument, makes his own figure or voice seem similar to yours, that kind of fantastic art is called mimetic.'[19]

Thus, in mimetic art the instrument of imitation is the imitator himself, while in the nameless opposite art the instrument of imitation is other than the imitator.
In the seventh division-step, described in 268a9-c4, the dissembler is divided into the demagogue and the sophist, the demagogue being
'one who can dissemble in long speeches in public before a multitude',
while the sophist is someone
'who does it in private in short speeches and forces the person who converses with him to contradict himself.'
Thus, he who listens to a dissembler is deceived, if that dissembler is a demagogue, and contradicted not by himself but by another instrument of deceit (namely the demagogue himself), while if the dissembler is a sophist, the listener is forced by the sophist to himself become the instrument of deceit.
Thus we have:

*(f4) Proposition.* The ratio of the fourth step in the Division of the Sophist, namely the ratio of the imitator who uses other instruments of imitation to the imitator who is himself the instrument of imitation
is equal to the ratio of the seventh step in the Division of the Sophist, namely the ratio of the demagogue, whose listener is contradicted and deceived not by himself but by another, to the sophist, whose listener is forced to be contradicted and deceived by himself.
Propositions (f2) and (f4) provide powerful additional evidence in favor of the anthyphairetic interpretation of Division and Collection.

## 5. A PLATONIC BEING IS AN ANTHYPHAIRETIC SELF-SIMILAR ONE

We are now ready to describe the precise sense in which a Platonic Form or Being can be regarded as a One.

*(a) What kind of a One for the Platonic Being?*
Plato distinguishes two states of Beings in the *Sophistes* 255c-d: the more exalted 'Being itself' and the more lowly Being 'with respect to' ('pros ti'), a relative Being. Only the latter is accessible and knowable to the human intellect, and we are here concerned exclusively with this Being. It is clear that, for Plato, the principal property of a Platonic Being—variously described as a One or a Unity, as simple, partless or indivisible—is precisely its Oneness. But what kind of Oneness?

*(b) The kinds of One that are rejected for a Platonic Being*

*(b1) A Platonic Being is not the absolutely partless One.* The Absolute One, the One without any parts whatsoever, the really partless and indivisible One which corresponds to the One of the first hypothesis (137c-142a) in the *Parmenides*, is explicitly rejected as Being (especially in 141e3-142a1). In this rejection, Plato follows Zeno, *Fragment B2*.

*(b2) A Platonic Being is not the cumulative One.* Cornford[22] thought that Collection is just the inverse of Division: he believed we could obtain the initial Genus by summing all the pieces of a Division together. But there are serious problems with this interpretation of Collection, and the corresponding interpretation of a Platonic Being as a One that is simply the sum of an infinity of parts, a summed into a totality. To begin with, there is no periodicity involved in such a Collection, although, as we note in 5 (d2) below, every Platonic description of Collection stresses its periodic nature. Secondly, the type of Collection suggested by Cornford is explicitly rejected both in the *Parmenides* 128e-130a (especially 129c4-d6 given in 5(c), below) and in the *Philebus* 14d-e, as describing sensible and not intelligible entities:

'and this sort of thing also should be disregarded, when a man in his talk dividing the members, and at the same time the parts, of anything, acknowledges that they all collectively are that one thing'[23] *(Philebus* 14e).

---

[22] F. M. Cornford, «Mathematics and Dialectic in the Republic VI-VII», *Mind N.S.*, 41, No. 161 (Jan. 1932), p.37-52 and No. 162 (Apr. 1932), p.173-190.; F. M. Cornford, *Plato's Theory of Knowledge, the Theaetetus and Sophistes of Plato*, Routledge & Kegan Paul, London, 1935.; F. M. Cornford, *Plato and Parmenides*, Routledge & Kegan Paul, London, 1939.

*(c) A Platonic Being is described as a self-similar One in the Sophistes and in the second hypothesis of the Parmenides.* At first sight, infinite divisibility seems to make a mockery of—and to run against—any decent notion of One, as it produces an infinite number of parts, namely the remainders at each stage of the division process. But as it is made clear in the *Sophistes* 244d-245b, 257c, 258e, the *Parmenides* 129c-d, and the *Parmenides* Second Hypothesis, divisibility is no obstacle for Oneness, specifically for self-similar Oneness in which every part is the same as the whole:

> 'But yet nothing hinders ('ouden apokoluei') that which has parts ('to memerismenon') from possessing the attribute of unity ('to pathos …tou henos') in all its parts ('en tois meresin pasin'), and in this way every Being and Whole ('pan te on kai holon') to be One ('hen')'[19] (*Sophistes* 245a1-3);

> 'but what is surprising if someone shall show that I am one and many ('hen kai polla')? When he wishes to show that I am many, he says that my right side is one thing and my left another, that my front is different from my back, and my upper body in like manner different from my lower; for I suppose I have a share of multitude. To show that I am one, he'll say I am one man among the seven of us, since I also have a share of the one. So he shows both are true. Now, if someone should undertake to show that sticks and stones and things like that are many, and the same things one, we'll grant that he has proved that something is many and one ('polla kai hen'), but not that the one is many ('to hen polla') or the many one ('ta polla hen'): he has said nothing out of the ordinary, but a thing on which we all agree.' (*Parmenides* 129c4-d6);

> 'there are infinitely many parts of Being' (*Parmenides* 144b6-7), 'to each of the parts of Being befits the One' (*ibid.* 144c6), 'neither the being is lacking in relation to the one, nor the one in relation to the being, but they are equalized being two for ever in all ways'[24] (*ibid.* 144e1-3).

*(d) A Platonic Being is an Anthyphairetic Self-similar One*

*(d1) Anthyphairetic Division of a Platonic Being into an infinite multitude of parts.* This relative Being possesses, in analogy to a mathematical ratio of magnitudes, two opposite to each other parts called powers in the

---

*Theaetetus* 156a6-b1—$B_1$ and $A_1$—one of which—say $A_1$—acting and the other—$B_1$—being acted on (247d-e); this is the Platonic entity, called by Aristotle an 'indefinite dyad' in, for example, the *Metaphysics* 1082a13-15, 1083b35-36.

Examples of Platonic 'indefinite dyads' defining a Being are: {$A_1$ beautiful, $B_1$ not beautiful}(257d-e); {$A_1$ great, $B_1$ not great} and {$A_1$ just, $B_1$ not just}(258a); {being $A_1$, non-being $B_1$} (258b-c); {$A_1$ self-restraint and $B_1$ bravery}, the indefinite dyad of the Form Virtue (mentioned and analysed in the final section of the *Politicus* 305e-311c); and {$B_1$ One, $A_1$ Being}, the indefinite dyad of the Form One Being (in the second hypothesis of the *Parmenides* 142b1-143a3). Indefinite dyads without an explicit reference to a Platonic Being, such as the sensible ones {cold, warm}, {fast, slow}, {more, less}, appear in the *Philebus* 23b5-25e3 and are instances of the Infinite ('apeiron'), interpreted in Negrepontis[25] as a philosophic form of infinite anthyphairesis. The initial indefinite dyad of the Beings the Angler and the Sophist is, in each case, the pair of the opposite species-parts {$B_1$, $A_1$} produced at the first stage by dividing the initial genus G.

As shown in the second hypothesis of the *Parmenides* 142b1-143a3, for the dyad {One, Being}, the initial infinite dyad {$B_1$, $A_1$} produces, by anthyphairetic division, the infinite sequence of parts-species $B_1, A_1, B_2, A_2, \ldots, B_k, A_k, B_{k+1}, A_{k+1}, \ldots$. In the abbreviated form of the Division, as it appears in the Divisions in the *Sophistes*, the parts of a Platonic Being are produced by dividing the initial Genus G at the first stage into two species B1 and A1, one of which—the one which is a true predicate to the Platonic Being, $A_1$—turns into a genus and is divided again at the second stage into two opposite species-parts, say $B_2$ and $A_2$. Because of Logos-periodicity, the Division continues ad infinitum, producing the infinite sequence of parts-species. The multitude of the parts of a Platonic Being is thus infinite.

*(d2) Anthyphairetic Periodicity and the finite number of Logoi of a Platonic Being.* A Logos of a Platonic Being is the ratio of two successive species-parts, namely either the ratio of the opposite species at some level, say k, $B_k/A_k$, or a ratio of the form $A_k/B_{k+1}$. By the Logos Criterion and the resulting periodicity, there are natural numbers m>n, such that $B_n/A_n = B_m/A_m$. Thus, the multitude of Logoi in a Platonic Being is finite and consists of all the Logoi in a period, namely:
$B_n/A_n, A_n/B_{n+1}, B_{n+1}/A_{n+1}, \ldots, B_{m-1}/A_{m-1}, A_{m-1}/B_m$ (and $B_m/A_m = B_n/A_n$).

The finiteness of the number of Logoi is clearly confirmed in 257a4-6. It is true that Collection is not associated in a crystal-clear manner with periodicity in the *Sophistes*; but this association is definitely clarified elsewhere in Plato, for example in the *Politicus* 283b-287b, in the *Philebus* 14c-19a, and principally in the *Phaedrus* 264e-266c.

---
[25] S. Negrepontis, «The Anthyphairetic nature of the Platonic Principle....»...

*(d3) Participation of a part in a Logos.* Each part-species participates in a Logos. Thus the species Bk participates in the Logos $B_k/A_k$ and the species Ak participates in the Logos $A_k/B_{k+1}$ or in the Logos $B_k/A_k$. This is precisely the basis of Platonic participation of the Sensibles (closely connected with the parts of the Division) to the Intelligibles (essentially the Logoi of Collection).

*(d4) Collection reveals the Oneness of the Platonic Being.* The Divisions-and-Collections described in Sections 3 and 4, which occur respectively at the start (Angler) and the end (Sophist) of the *Sophistes* dialogue, are meant to be examples of the 'name and logos' method which, according to our analysis, we have interpreted to mean 'division and philosophic Logos Criterion', namely 'philosophic anthyphairetic expansion and Logos Criterion'. In the long stretch between the two Divisions (*Sophistes* 236d-260b), the emphasis shifts from the method of 'name and logos' to the method of 'division and collection' or, equivalently, the 'division according to kinds' described in 252b1-5 and 253d1-e3. Since 'name' is clearly associated with division, the new element in this description is 'collection' ('sunagoge'), and it is quite clear that it must be connected with the 'Logos' in the previous description, 'name and logos'. Of course, 'division' is 'division of the one into many' and 'collection' is 'collection of the many into One', so 'Logos' must somehow have the power to collect the many into One.

*(d5) Similarity of the Logoi.* By the philosophic analogue of the definition of the proportion of magnitudes, given in the *Topica,* and the resulting Logos Criterion, the set of Logoi in the anthyphairesis of each Logos k—say $B_k/A_k$ or $A_k/B_{k+1}$—is the same, except that it appears in a cyclic permutation of the original (cf. 1(i) above). In this way, we regard all the Logoi of a Platonic Being as the same and similar to each other.

*(d6) Equalization of all parts of a Platonic Being.* The fundamental property of a Platonic Being is that the two opposite parts—being a and non-being b—of the indefinite dyad possess the 'power' (247d-248d), which is clarified as the 'power to communicate' with each other (252d3, 253a8, 253e1, 254b8, 254c5) and still further clarified to be a 'power for equalisation with each other' (257b). The same basic property of a Being, is described in the *Politicus* 305e-311c with regard to the Being or Form of Virtue, where the dyad **a** bravery and **b** self-restraint is called the Logos that defines Virtue.

The way the two opposing parts of the indefinite dyad of a Being (and in fact any two successive remainders in the sequence a, b, $a_1$, $b_1$,..., $a_k$, $b_k$, $a_{k+1}$,... of the anthyphairetic division) become equalised and collected into One by Logos is well-explained in the *Sophistes* 257b1-260b3, especially in 257b, 257e9-258c5: it is by considering not the opposing parts themselves (which are 'not-beings', 256d-e, 257a), but the corresponding Logoi-ratios (which are relative Beings) into which they participate; thus, we consider instead of the non-being a, the relative Being the logos a/b, and instead of the non-being b, the relative Being logos $b/a_1$ (257e2-258b5).

The interpretation of 'Collection' as precisely this 'equalisation of parts' becomes clear in 253d1-e3, where the method of 'division into kinds', which is clearly equivalent to 'division and collection', is declared equivalent to 'the power of the parts to communicate'—namely, as we have already seen, to 'the power of the parts to be equalised with each other'.

Thus, in a Relative Being or Logos, there is an infinite number of parts or 'not-beings' (as confirmed in 256e5-7), but all parts are equalised in the sense of fractal self-similarity.

*(d7) Note.* The interpretation of the *Sophistes* presented in this paper is completely novel. Most of the received interpretations, such as those given by Cherniss[26], Cornford[27], Ackrill[28], Taylor[29], Owen[30], Vlastos[31], Frede[32] and Dorter[33], to name some of the most prominent, have very little in common with our mathematical analysis; specifically, there is no sign or suspicion in their interpretations that periodic anthyphairetic expansions might have a role to play, or that a non-being is just a part in the division which is turned into a fractal self-similar being through its participation in the logos, formed with the immediately succeeding part, their attention having been sidetracked away from mathematics and mostly to the syntax of the verb to be.

---

[26]H.F. Cherniss, *The Philosophical Economy of the Theory of Ideas*, 1936, reprinted in R.E. Allen (ed.), *Studies in Plato's Metaphysics*, Routledge & Kegan Paul, London, 1965; H. Cherniss, *The Riddle of the Early Academy*, University of California Press, Berkeley, 1945, p. 53-55; H. Cherniss, *Plato as Mathematician*,.....

[27] F.M. Cornford, «Mathematics and Dialectic....»; F.M. Cornford, *Plato's Theory of* .....; F.M. Cornford, *Plato and Parmenides*…

[28]J.L. Ackrill, «SUMPLOKE EIDON», (1955), reprinted in R.E. Allen (ed.), *Studies in Plato's Metaphysics*, Routledge & Kegan Paul, London, 1965, p. 199-206; J.L. Ackrill, «Plato and the Copula: Sophist 251-259», (1957), reprinted in R.E. Allen (ed.), *Studies in Plato's Metaphysics*, Routledge & Kegan Paul, London, 1965, p. 207-218.

[29] A.E. Taylor, *Plato: the Man and his Work*, Methuen, London, 1960.

[30] G.E.L. Owen, «*Plato on Not-Being*», in G. Vlastos (ed.) *Plato I: Metaphysics and Epistemology*, Anchor Books, Garden City, 1970, p. 223-267.

[31] G. Vlastos, «An Ambiguity in the Sophist», in *Platonic Studies*, 1970, Second Edition, Princeton University Press, Princeton, 1981, p. 270-322.

[32] M. Frede, «Plato's Sophist on false statements», reprinted in R. Kraut (ed.), *The Cambridge Companion to Plato*, Cambridge University Press, Cambridge, 1992, p. 397-424.

[33] K. Dorter, *Form and Good in Plato's Eleatic Dialogues: the Parmenides, Theaetetus, Sophist, and Statesman*, University of California Press, Berkeley, 1994.

Palmer[34] develops an interpretation of the Second Hypothesis of the *Parmenides*, according to which the infinite divisibility of the One Being is not an obstacle, but rather an asset, to its nature as a One, thanks to 'the demonstration of the one being's resistance to the type of division envisaged here via a demonstration of how its predicate parts "one" and "being" persist through every possible division' (p. 225). He thus approaches the self-similar Oneness of a Platonic Being, but in an incomplete way, taking into account only the infinite divisibility (which cannot by itself produce self-similarity) and not the anthyphairetic periodicity.

## 6. THE ANTHYPHAIRETIC INTERPRETATION OF PLATONIC NUMBERS

The self-similarity of a Platonic Being opens up the way for the definition of Platonic numbers: the anthyphairetic interpretation of Platonic numbers is a natural by-product of the anthyphairetic interpretation of the method of Division and Collection. We shall briefly describe the nature of the Numbers in accordance with the anthyphairetic interpretation.

These are certainly the same numbers described in *Philebus* 56c10-e6 (where they are differentiated from the numbers used by the many) and the *Politeia* 522c1-526c7 (where they are differentiated from the numbers as used by, say, Agamemnon). The differentiation from common usage is twofold: a Platonic number consists of units that are (i) entirely equal to each other, and (ii) divisible into an infinite multitude of parts.

The nature of these numbers is well described by Aristotle in *Metaphysics* 987b25-988a1, where it is stated that:

(a) 'the number two consists of the two elements of an indefinite dyad, equalized ('isasthenton') by the principle of the One' (1081a23-25, 1083b30-32, 1091a24-25),

(b) contrary to the mathematical units, the units of the Platonic ('eidetic') numbers are equal for a fixed number, but different for different eidetic numbers (1080a23-30), and

(c) the species are numbers (987b20-22).

The difficulties modern Platonists have encountered in reconstructing the Platonic numbers according to Aristotle's requirements has led to a downgrading of Aristotle's account of Platonic numbers (misunderstanding on his part, unwritten dogmas, and so on). But Plato

---

[34] J. Palmer, *Plato's reception of Parmenides*, Clarendon Press, Oxford, 1999, Chapter 10.

describes Platonic numbers in the second hypothesis of the *Parmenides*, and the anthyphairetic interpretation of this hypothesis, discovered by S. Negrepontis[35], allows Aristotle's description of Platonic numbers to be confirmed. Thus, Plato's definition of the number two in the *Parmenides* 144b-e as the pair of One and Being fully agrees with (a), since the self-similar Oneness of the Platonic Being precisely assures that these two successive units-species are equalized ('exisousthon', 144e2) by periodicity and Logos (or, as Aristotle calls it, the principle of the One). The units for forming a Platonic number are the species, numbered according to their successive generation in the anthyphairetic Division of the Being, precisely as in the *Parmenides* 144b-e, the *Philebus* 16c5-17a5 and the *Sophistes* 258c3, where each not-being part-species is claimed to be 'one species among the many Beings possessing number ('enarithmon'). Platonic numbers thus indeed are species, as Aristotle claims (c). The numbers in the Being Sophist are the numbered species in

Table 10. The Platonic Numbers in the Platonic Being Sophist

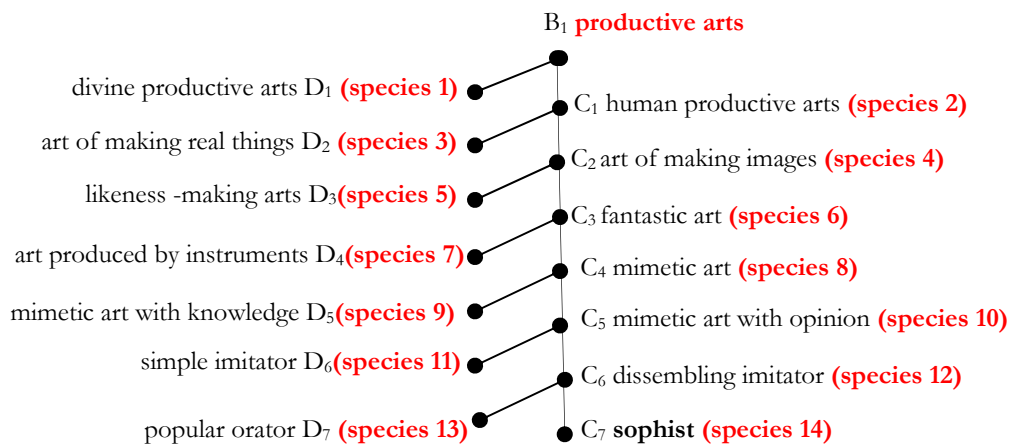

According to this interpretation, Platonic numbers are not absolute but only relative to a Platonic Being: the initial unit (hence (b)). In fact, however, according to the *Parmenides* 148d5-149d7 and the *Sophistes* 257a4-6 ('hosaper…tosauta', like the 'tosauta hosaper' in the *Parmenides* 144d5, both following the 'tosauta…hosa' in Zeno's *Fragment B3*), the number of units in any Platonic Being is finite and equal to (the number of different Logoi)+1 (this number being determined by the length of the period).

## 7. THE ANTHYPHAIRETIC INTERPRETATION OF TRUE JUDGEMENT

Once we have obtained the anthyphairetic interpretation of Division and

---

[35] S. Negrepontis, «The periodic Anthyphairetic….»….

Collection, it is relatively straightforward to realize that, for Plato, True Judgement is an anthyphairetic approximation of the Division of a Platonic Being. We will outline the powerful evidence for such an interpretation.

*(a) Knowledge of an intelligible Being is True Judgement and Logos.* Knowledge of intelligible Beings is described by Plato as 'True Judgement with logos' in the *Sumposium* 202a2-10, 209e5-211d1, in *Meno* 97e2-98c4 (where instead of 'logos', the equivalent expression 'logismos aitias' is employed), and in the *Theaetetus* 201c8-d3, 202b8-d7. The same knowledge of intelligible Beings is also described equivalently by Plato as 'Name and Logos' in the *Theaetetus* and the *Sophistes*, and as 'Division and Collection' in the *Sophistes*. It is further clear that the two descriptions 'Name and Logos' and 'Division and Collection' are equivalent in the sense that 'Name' and 'Division' have the same meaning, while 'Logos' and 'Collection' are also synonymous.

*(b) True Judgement is between Knowledge and Ignorance and it is alogos; hence True Judgement is related only to Division.* According to the *Sumposium* 202a2-10, True Judgement is not to be equated with Ignorance, since it is a part of the knowledge of a Being, but is also not Knowledge, because it is an 'alogos' entity and 'has no power to provide logos'. Thus, True Judgement lies between Knowledge and Ignorance. A similar description of True Judgement occurs in the long passage *Politeia* 477e9-479e6, where it is stated that True Judgement is 'neither Knowledge nor Ignorance' but something 'darker than Knowledge and with more light than Ignorance', meaning that it is located 'between' Ignorance and Knowledge.

The precise reason why True Judgement is something less than Knowledge lies in the fact that True Judgement is 'alogos'. This is stated clearly in the *Symposium*, as we have seen, but also in *Theaetetus* 201c8-d3 and in *Timaeus* 28a, 51d3-e4.

But, as we have seen, Knowledge of intelligible Beings is equivalently described as Name-Division and Logos-Collection. Since True Judgement is 'alogos', it follows that True Judgement is related ONLY to Division.

*(c) Judgement is a finite process.* On account of (b), the most natural conclusion would be to equate True Judgement with Division. However, Division is infinite, while the two descriptions of Judgement in Plato clearly show that it is a finite process. Thus, in *Theaetetus* 189e4-190a6: 'But the soul, as the image presents itself to me, when it thinks, is merely conversing with itself, asking itself questions and answering, affirming and denying. When it has arrived at a decision ('epaxasa'), whether slowly or with a sudden bound, and is at last agreed, and is not in need of a dichotomical decision ('me distaze'), we call that its Judgement'; and in *Sophistes* 264a9-b1: 'thought is conversation of the soul with itself, and Judgement is the final result ('apoteleutesis') of thought'.

There is considerable backing for interpreting 'distaze' as a dichotomical decision:

(i) Hesuchios *Lexicon,* in Delta 1976 and 2021, [distazei] is rendered as 'dichonoei' (being of a mind divided in half), while in Delta 2086, 'distazei' is equated with 'amphidoxei' (having both opinions);

(ii) *Magnum* 278,35-36 explains 'distazein' as thinking in a double manner ('dichos bouleuesthai') or as 'to have two stations in thinking ('duo staseis echein kata dianoian');

(iii) Philoponus notes in *eis Aristotle's On Soul* 15,78,11-26 that Judgement can be compared to 'distazei', which is likened to travelling a divided road ('schisten odon'), considering to go in one or in the other direction.

'Me distazei' thus refers to a condition in which there are no more divisions, meaning the process—which, as we saw in 5(b), is a process of division—has come to an end.

*(d) True Judgement has an approximative character, comes in degrees, and in succession.* The description of Judgement in the *Cratulus* 420b6-d2 as an arrow directed in pursuit of its target makes clear its approximative nature; as Judgement gets nearer to its target, it achieves better approximations.

In the *Politeia* 506c6-10, even 'the very best Judgements ('hai beltistai') are called blind (no doubt because they are 'alogoi'); a hierarchy of Judgements from bad to better is clearly implied.

In the *Theaetetus*, True Judgement is described, in 191c8-d2, in terms of a block of wax in our souls whose quality, varying from soul to soul, is a basic factor in the quality of the Judgement. Thus, the wax can be larger ('meizon') or smaller ('elatton'), purer ('katharoterou') or more impure ('koprodesterou'), harder ('skleroterou') or softer ('hugroterou') or, finally, of proper quality ('metrios echontos'). This is further detailed in 194c4-195b1: the wax may be deep ('bathus'), abundant ('polus'), smooth ('leios') or properly kneaded ('metrios orgasmenos'), in which case the imprints will be clear ('kathara', 'saphe'), being of sufficient depth ('hikanos tou bathous echonta') longevity ('poluchronia') or spaciousness ('en euruchoria'); in contrast, the wax can also be shaggy ('lasion'), unclean ('koprodes'), impure ('me katharou'), overly soft ('hugron sphodra') or hard ('skleron'), rough ('trachu'), stony ('lithodes'), mixed with earth or dung ('ges e koprou summigeises'); correspondingly, the imprints on it may become indistinct ('asaphe', 'eti asaphestera'), either by lacking depth ('bathos ouk eni'), melting together ('sugcheisthai') quickly becoming blurred ('amudra'), or by being piled one on top of another due to lack of space ('ep' allelon sumpeptokota…hupo stenochorias').

Another strong indication of the approximative nature of True Judgement is found in the dialectical kernel of Diotima's lecture on the love of beauty in the *Symposium* 209e5-211d1; Diotima's lecture is a vivid description of the ascent to a unitary Platonic Being in successive ('ephexes') approximating levels of True Judgement, crowned with the sudden ('exaiphnes') appearance of Logos.

> 'If you pursue with right judgement ('orthos metie') the lesser mysteries of love, they may lead to the perfect and collecting ones ('telea kai epoptika')' (209e5-210a2); 'He who would proceed with right judgement ('ton orthos ionta') towards this entity, should begin in youth to seek the company of corporeal beauty; and first, if his instructor has right judgement ('orthos hegetai'), to love one beautiful body and generate beautiful 'logoi' ('gennan logous kalous')' (210a4-8); and then, if the form of beauty ('to ep' eidei kalon') is his pursuit ('diokein'), how foolish would he be not to judge ('hegesthai') that the beauty in every body is one and the same. And when he perceives this…he will despise and judge ('hegesamenon') a small thing, and will became a steadfast lover of all beautiful bodies.' (210a8-b6); 'In the next stage he will judge ('hegesasthai') that the beauty of the soul is more precious ('timioteron') than the beauty of the outward form' (210b6-7); 'until he is compelled next to contemplate and see the beauty in institutions and laws and to understand that the beauty of them all is of one family, so that he will judge ('hegesetai') that corporeal beauty is a trifle' (210c3-6); 'and after institutions his guide will lead ('agagein') him on to the sciences' (210c6-7); 'He who has been instructed thus far in the things of love, and who has learned to see the beautiful in succession and with right judgement ('ephexes te kai orthos'), when he comes ('ion') toward the end will suddenly perceive a nature of wondrous beauty' (210e2-5), which is a Platonic Being, possessing unity and simplicity ('everlasting being' ('aei on') [211a1], 'simple everlasting being' ('monoeides aei on') [211b1-2], 'simple' ('monoeides')[36] (211e4)).

In 211b5-d1, in recapitulation, the process of approximation is called right judgement in love ('orthos paiderastein'), going to the things of love with right judgement ('orthos epi ta erotica ienai'), and employing the various levels of True Judgement as steps ('epanabasmois') (cf. the description of hypotheses as 'epibaseis kai hormas' toward the 'anhupotheton in the *Politeia* 511b6).

---

[36] Based on Plato, *The Symposium*, translated by Benjamin Jowett, Sphere Books, London, 1970.

*(e) [True] Judgement is the philosophic analogue of Finite Anthyphairesis.* Anonymous *Scholion* X.2 to Euclid's *Elements* describes Judgement as the philosophic analogue of Finite Anthyphairesis. Indeed, *Scholion* X.2 first states that it is in the nature of numbers to have finite anthyphairesis: 'those [the numbers] when divided they end in a definite end, unity' (cf. Propositions VII.1,2 of the *Elements*); this is so, the *Scholion* continues, because numbers act in a Judgement-like manner ('tous de arithmous doxastikos'), and Judgement is of a finite nature ('peperastai gar mallon he doxa').

*(f) True Judgement of a Platonic Being is a finite anthyphairetic approximation of the Division of that Being, analogous to the side and diameter numbers approximating the diameter to the side of a square.* As noted in (b), the fact that True Judgement and Logos is a description equivalent to Name and Logos and to the Division and Collection of Knowledge of a Platonic Being, coupled with the description of True Judgement as totally devoid of Logos ('alogos'), leads to the conclusion that True Judgement is related only to Platonic Division. Our interpretation of Platonic Division as the philosophic analogue of infinite anthyphairesis (in a Platonic Being), given in Sections 2 and 3, leads to the conclusion that True Judgement is related to the anthyphairetic Platonic Division of a Platonic Being.

Our first impulse would be to identify True Judgement with Platonic Division, but the evidence presented in (c) precludes such an identification: True Judgement is finitistic in character; more precisely, it is a process of finite anthyphairetic-like division. This conclusion is strengthened by the *Anonymous Scholion* X.2 to Euclid's Elements given in (e). Finally, it is clear from the copious evidence presented in (d) that True Judgement has an approximative nature and comes in degrees of higher and lower approximation (of the Platonic Being).

At this point, we will employ our knowledge of the approximation of ratios in ancient Greek Mathematics. It is well-documented in Fowler[37] that ALL known approximations of a ratio by a simpler one are anthyphairetic: namely, the anthyphairesis of the simpler ratio is an initial segment of the anthyphairesis of the approximated ratio. The Pythagoreans had perfect knowledge of the 'convergents'—the 'side and diameter numbers' as they called them—of the anthyphairesis of the diameter to the side of a square. Since True Judgement is an approximation of the Division of a Platonic Being, and since Platonic Division is anthyphairetic in nature, it is clear that True Judgement is precisely the philosophic analogue of the (finite) anthyphairetic approximation of the (infinite) anthyphairetic Platonic Division of a Platonic Being. Thus, the infinite approximative hierarchy of True Judgements is the analogue of the sequence of side and diameter numbers (in agreement with (e))[38].

---

[37] D. Fowler, *The Mathematics of…*, Chapter 2.

Table 11. True Judgement of level k=5 about the Platonic Being Sophist

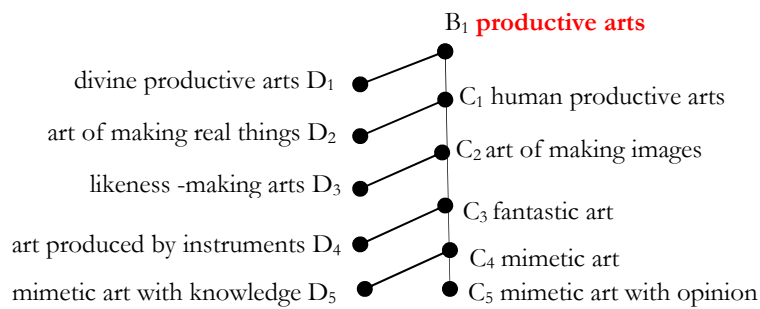

This true judgement is expressed in common language with the description 'Sophistry is a mimetic art with opinion'.

Table 12. True Judgement of level k=6 about the Platonic being Sophist

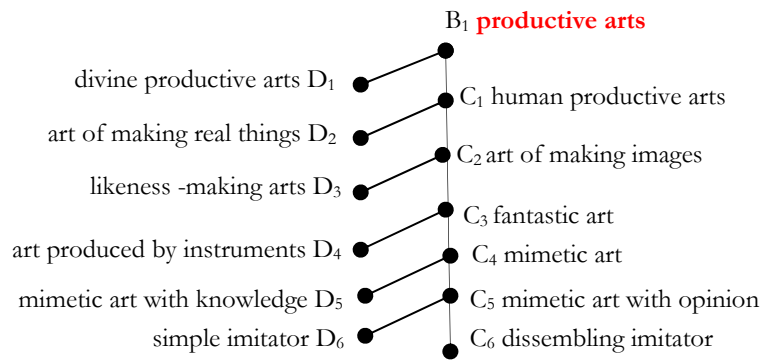

This true judgement is expressed in common language as the description 'Sophistry is an art of dissembling imitation'. The statement 'Theaetetus is sitting', the example of True Judgement in the *Sophistes* 263a2, is precisely such a True Judgement (about Theaetetus).

The interpretation of True Judgement, given here, will be employed for the interpretation of Plato's Theory of Falsehood, in the joint work[39].

---

[42]These are the 'natural logoi' ('phusikoi logoi') to which Proclus refers in *eis Parmeniden* 877, 32-885, 32 in relation to True Judgement (and the Third Man Argument).
[39] S. Negrepontis and S. Birba-Pappa, «Plato's Theory….»…

*Department of Mathematics,*
*Athens University, Athens 157 84, Greece*
snegrep@math.uoa.gr